\documentclass[11pt,leqno]{article}
\usepackage{graphicx, amsfonts, amsthm, amsxtra, verbatim, multicol}
\usepackage[mathscr]{euscript}
\begin{document}

\begin{center}
{\LARGE\bf Commuting holonomies and rigidity of holomorphic foliations 
}
\footnote{
Math. classification: 57R30, 14D99, 32G34 
} \\

\vspace{.25in} {\large {\sc 
Hossein Movasati}}\footnote{The first author is supported by the Japan Society for the 
Promotion
of  Sciences.},
{\large {\sc Isao Nakai}} \\
\begin{abstract}

In this article we study deformations of a holomorphic foliation
with a generic non-rational first integral in the complex plane. We consider two vanishing cycles 
in a regular fiber of the first integral  with a non-zero self intersection and with vanishing paths which intersect 
each other only at their start points.  It is proved  that 
if the deformed holonomies 
of such vanishing cycles commute then the deformed foliation has also a first integral.
Our result generalizes a similar result of Ilyashenko on the rigidity of holomorphic foliations with
a persistent center singularity.  The main tools of the proof 
are Picard-Lefschetz theory and the theory of iterated integrals for such deformations.


\end{abstract}
\end{center}
\newtheorem{theo}{Theorem}
\newtheorem{lemma}{Lemma}
\newtheorem{exam}{Example}
\newtheorem{coro}{Corollary}
\newtheorem{defi}{Definition}
\newtheorem{axio}{I}

\newtheorem{prob}{Problem}
\newtheorem{lemm}{Lemma}
\newtheorem{prop}{Proposition}
\newtheorem{rem}{Remark}
\newtheorem{conj}{Conjecture}
\def\End{{\mathrm End}}              
\def\hol{{\mathrm Hol}}        
\def\sing{{\mathrm Sing}}            
\def\cha{{\mathrm char}}             
\def\Gal{{\mathrm Gal}}              
\def\jacob{{\mathrm jacob}}          
\newcommand\Pn[1]{\mathbb{P}^{#1}}   
\def\Z{\mathbb{Z}}                 
\def\ZZ{\mathbb{Z}}               
\def\Q{\mathbb{Q}}                   
\def\C{\mathbb{C}}                   
\def\R{\mathbb{R}}                   
\def\N{\mathbb{N}}                   
\def\Fi{U}                   

\def\A{\mathbb{A}}                   
\def\uhp{{\mathbb H}}                
\newcommand\ep[1]{e^{\frac{2\pi i}{#1}}}
\def\Mat{{\mathrm Mat}}              
\newcommand{\mat}[4]{
     \begin{pmatrix}
            #1 & #2 \\
            #3 & #4
       \end{pmatrix}
    }                                
\newcommand{\matt}[2]{
     \begin{pmatrix}                 
            #1   \\
            #2
       \end{pmatrix}
    }
\def\ker{{\mathrm ker}}              
\def\cl{{\mathrm cl}}                
\def\dR{{\mathrm dR}}                

\def\hc{{\mathsf H}}                 
\def\Hb{{\cal H}}                    
\def\GL{{\mathrm GL}}                
\def\pedo{{\cal P}}                  
\def\PP{\tilde{\cal P}}              
\def\cm {{\cal C}}                   
\def\K{{\mathbb K}}                  
\def\F{{\cal F}}                     
\def\M{{\cal M}}
\def\RR{{\cal R}}
\newcommand\Hi[1]{\mathbb{P}^{#1}_\infty}
\def\pt{\mathbb{C}[t]}               
\def\W{{\cal W}}                     
\def\Af{{\cal A}}                    
\def\gr{{\mathrm Gr}}                
\def\Im{{\mathrm Im}}                
\newcommand\SL[2]{{\mathrm SL}(#1, #2)}    
\newcommand\PSL[2]{{\mathrm PSL}(#1, #2)}  
\def\Res{{\mathrm Res}}              

\def\L{{\cal L}}                     
\def\Aut{{\mathrm Aut}}              
\def\any{R}                          
\newcommand\ovl[1]{\overline{#1}}    
\def\T{{\cal T }}                    
\def\tr{{\mathsf t}}                 
\def\per {{\sf pm}}
\def\ring{{\sf R}}
\def\BM{{\sf H}}
\def\OO{{\cal O}}
\def\rank{{\rm rank }}
\section{Introduction}
In a deformation of an integrable foliation one obtains the first Melnikov function as 
an Abelian integral whose zeros give rise to limit cycles  
in the deformed foliation, see for instance \cite{mov0}. In the case where the
Abelian integral is  identically zero such limit cycles are controlled 
by higher order Melnikov functions  and  L. Gavrilov in \cite{ga06} has shown 
 that such functions
can be expressed in terms of iterated integrals and so they
 satisfy certain Picard-Fuchs equations. 
 In a different context, the second named author and K. Yanai in \cite{nakai} 
 have used iterated integrals
 to investigate the existence of relations between formal diffeomorphisms.
 Basic properties 
 of iterated  integrals were established by A. N. Parsin in 1969 and a systematic approach
 for de Rham cohomology type theorems for iterated integrals 
 was made by K.-T. Chen around 1977. In the present text we use iterated integrals and 
investigate the non-existence of non-trivial commutator 
relations between deformed holonomies.
 
 Let us consider a polynomial in two variables  $f(x,y)\in \C[x,y], \deg(f)\geq 3$ and perform a 
 perturbation
\begin{equation}
\F_\epsilon: df+\epsilon \omega,\ \epsilon\in (\C,0),\ \deg(\omega)\leq \deg(f)-1,
\end{equation}
where $\omega=Pdx+Qdy$ is a polynomial differential form, $(\C,0)$ is a small neighborhood of the origin in $\C$ and $\deg(\omega)$ 
 is the maximum of $\deg(P)$ and $\deg(Q)$.
We take a path $\delta\in \pi_1(f^{-1}(b),p)$, where $b$ is a regular value of $f$ and $p$ is a 
point in $f^{-1}(b)$, and  ask for the conditions on $\omega$ such that 
the deformed holonomy 
$h_\epsilon:\Sigma\rightarrow \Sigma$, where $\Sigma$ is a transversal section to 
$\F_0$ at $p$,  is identity.
The first result in this direction is due to Yu. Ilyashenko:

Assume that the last homogeneous piece of $f$ is a product of $\deg(f)$ distinct lines
and the critical points of $f$ are non-degenerate with distinct images. These are generic
conditions on $f$.
For simplicity assume that the coefficients of $f$ are real numbers and 
 take $\delta$ one of the oriented ovals which lies in the level curves of the map
$f:\R^2\to\R$ and call $\delta$ a vanishing cycle. 
\\
{\bf Theorem.}(Ilyashenko, \cite{il69})
{\it  
 If $h_\epsilon$ is the identity map for all 
$\epsilon\in(\C,0)$ and the homology class of 
$\delta$ is a vanishing cycle then there is a polynomial $g\in\C[x,y],\ \deg(g)\leq \deg(f)$ such that 
$\omega=dg$ and so $\F_\epsilon:d(f+\epsilon g)=0$ is again Hamiltonian.
}

A generalization of this result for pencils of type $\frac{F^p}{G^q}$ in $\Pn 2$ and pencils
in arbitrary projective manifolds and logarithmic foliations is done in the articles
\cite{mov0, mov}. The theory of iterated integrals gives us further generalizations of
the above theorem for cycles with zero homology classes. For $a,b$ in a group $G$
let $(a,b):=aba^{-1}b^{-1}$ be the commutator of $a$ and $b$.
\begin{theo}
\label{12.4}
For a generic polynomial $f$ as before, let us assume that $h_\epsilon$ is the identity map for all $\epsilon\in(\C,0)$, $\delta=(\delta_1,\delta_2)$ and
the homology classes of $\delta_1$ and $\delta_2$ vanish  along two paths which do not
intersect again  except at $b$. Further assume that the homology
classes of $\delta_1$ and $\delta_2 $ have non-zero intersection.
Then
$\omega$ is an exact form and so $\F_\epsilon$ is again Hamiltonian. 
\end{theo}
In \S \ref{second} we 
state Theorem \ref{commuting}  for strongly tame functions. Two special cases 
of  Theorem \ref{commuting} are Theorem \ref{12.4} and the following:
Let $M$ be a projective compact manifold of dimension two and $\F(\omega_0)$ 
be a holomorphic 
foliation in $M$ obtained by a generic 
non-rational Lefschetz pencil (see \cite{lam81} and \S \ref{first}). Here 
$\omega_0$ is a global holomorphic section of $\Omega_M^1\otimes L$ such that
the zero locus of $\omega_0$ is a finite set in $M$,
where $L$ is a line bundle on
$M$ and $\Omega_M^1$ is the sheaf of holomorphic $1$-forms in $M$. 
Let 
$$
\F_{\epsilon}=\F(\omega_0+\epsilon\omega_1),\epsilon\in (\C,0),\ \omega_1\in H^0(M,\Omega_M^1\otimes L)
$$ 
be a linear deformation of $\F(\omega_0)$.
If $\delta_1,\delta_2$ are two vanishing cycles with the same properties
as in Theorem \ref{12.4}, $H_1(M,\Q)=0$ and the holonomies associated to  $\delta_1$
and $\delta_2$ commute  then $\F_\epsilon$ has a first integral.

A typical example of the situation of Theorem \ref{12.4} is the following:  
Assume that $f:\R^2\rightarrow \R$ has two non-degenerate
critical points: $p_1$ a center singularity and $p_2$ a saddle singularity.
Assume that there is no more critical value of $f$  between $f(p_1)$ and $f(p_2)$.
The real vanishing cycle around $p_1$ and the complex vanishing cycle around $p_2$ satisfy
the hypothesis of Theorem \ref{12.4}. For a more explicit example take $f$ to be the product of $d$
degree $1$ real polynomials which are in general position and deform it in order 
to obtain a generic polynomial required by Theorem \ref{12.4}. For a precise description of
the generic properties we have posed on $f$ see \S \ref{first}.

The paper is organized as follows: In \S \ref{first}
we define a strongly tame function in an affine variety. 
In \S \ref{second} we state and prove Theorem \ref{commuting} which is a general form of Theorem \ref{12.4}.
\section{Deformation of Holomorphic foliations}
\label{first}
In this section we consider a smooth affine variety $U$ and 
a polynomial function $f$ in $U$ and look at it as a morphism $f:U\to \C$ of algebraic varieties.
There is a compact projective  manifold  $M$ and a divisor $D$ in $M$ such that 
$U=M\backslash D$. There is also
a rational morphism $\bar f:M\rightarrow \Pn 1 $ which coincides with $f:U\to \C$ when restricted to
$U$. For 
$t\in\C$ we define $U_t:=f^{-1}(t)$.
\begin{defi}\rm
The morphism $f:U\rightarrow \C$ is tame if
\begin{enumerate}
\item
The divisor at infinity  $D$ is smooth and connected and further we have 
$H_1(U,\Q)=0$;
\item 
The foliation $\F(df)$ is not rational, i.e. the closure of a generic fiber of $f$ is not isomorphic
to $\Pn 1$.
\item
$f$ has non-degenerate singularities
$p_i,\ i=1,2,3\ldots,\mu$ with distinct images $c_i:=f(p_i)$ (critical values of $f$);
\item
A generic fiber of  $f$ is connected and its closure in $M$ intersects $D$ transversally.
\end{enumerate}
\end{defi}
One usually calls $\mu$  the Milnor number of $f$ and 
$C:=\{c_1,c_2,c_3,\ldots,c_\mu\}$  the set of critical values of $f$.
Ehresmann's theorem implies that a tame morphism is topologically trivial over $\C\backslash C$. 
We have two main examples in mind. The first is 
a generic  Lefschetz pencil (see \cite{lam81}) in a projective manifold $M\subset \Pn m$.
The first and second conditions  become intrinsic properties of the pair $(M,\Pn m)$.
For this example, one can take $D$ in such a way  that $\bar  f$ is also topologically
trivial over $\infty$.
The second example is mainly used in planar differential equations.
Let $f$ be a polynomial in two variables with $\deg(f)=d\geq 3$. 
We may compactify $\C^2$ inside $\Pn 2$, and look at $\F(df)$ as a foliation in 
$\Pn 2$.  For a generic choice of the coefficients, the polynomial $f$ is tame. 
For instance, to obtain the fourth condition one assumes that 
$\{[x;y]\in \Pn 1_\infty \mid  f_d(x,y)=0\}$  has $d$ distinct points, 
where $f_d$ is the last homogeneous piece of $f$.
In this case  $D\cong\Pn 1$ is not a regular fiber of $f$. Geometrically seen, 
$d$ sheets of a regular fiber
of $f$ accumulate at $D$.  

We take a distinguished system of paths $\gamma_i,\ i=1,2,\ldots, \mu$ in $\C$ (see \cite{arn}).
The path $\gamma_i$ connects  a regular value $b$ of $f$  to 
$c_i$ and does not intersect
other paths except at $b$. Let $\delta_i\in H_1(U_b,\Z)$ be the vanishing cycle along $\gamma_i$.
One calls $\delta_i,\ i=1,2,\ldots,\mu$ a distinguished basis of vanishing cycles.
The Dynkin diagram of $f$ is  a graph whose vertexes are vanishing cycles $\delta_i$.
The vertex $\delta_i$ is connected to $\delta_j$ if and only if 
$\langle \delta_i, \delta_j \rangle\not =0$, where
$$
\langle \cdot,\cdot \rangle: H_1(U_b,\Z)\times H_1(U_b,\Z)\to\Z
$$
is the intersection form in $H_1(U_b,\Z)$.
The morphism $f$ is called strongly tame if $f$ is tame and its Dynkin diagram is connected.
 A generic Lefschetz pencil and a generic polynomial in two variables discussed above 
are strongly tame. For a proof see  \cite{lam81} 7.3.5 and \cite{hos2} Theorem 2.3.2, 2.
The polynomial case has been proved in \cite{il69}. It follows also from the following:
If a tame polynomial $f$ is obtained by a topologically trivial 
deformation of a morphism $g:U\rightarrow \C$ with only one singularity then
the Dynkin diagram of $f$ is connected and so it is strongly tame (see 
\cite{la74,ga74,gu77}).
By a topologically trivial deformation we mean the one in which the topological structure of
the smooth fiber does not change. 
 \begin{prop}
 \label{6mar06}
If $f$ is a strongly tame morphism then
\begin{enumerate}
\item 
A distinguished basis of vanishing cycles generate $H_1(U_b,\Q)$;
\item
For a cycle $\delta\in H_1(U_b,\Q)$ such that  $H_1(U_b,\Q)\to \Q,\ \delta'\mapsto \langle \delta,\delta'\rangle$ is not the zero map, the action of the monodromy  on 
$\delta$ generates $H_1(U_b,\Q)$. In particular, this is true for vanishing cycles. 
\end{enumerate}
\end{prop}
\begin{proof}
The first part can be proved by a slight modifications of the arguments of
\cite{lam81}, \S 5. For a precise proof see \cite{hos2} Theorem 2.2.1.
The second part follows from the first part, the connectivity of the Dynkin diagram and
 Picard-Lefschetz formula.
\end{proof}
Let $\Omega_\Fi^{i}$ be the set of meromorphic
differential $i$-forms in $M$ with poles along $D$.
A peculiar property of a tame polynomial is that if $\int_{\delta_t}\omega=0$ for a continuous
family of vanishing cycles $\delta_t$ and $\omega\in\Omega_U^1$ then $\omega$ is relatively
exact, i.e.
$
\int_{\delta}\omega=0,\ \forall \delta\in H_1(\Fi_t,\Z), t\in \C\backslash C
$
or equivalently $\omega$ is of the form $dP+Qdf$ for some $P,Q\in\Omega^0_U$
(see \cite{mov0} Theorem
5.1). 
 
 The Brieskorn module 
$$
H=\frac{\Omega^1_U}{df\wedge\Omega^0_U+d\Omega^0_U },
$$
is a $\C[t]$-module in a canonical way: $t[\omega]:=[f\omega]$.
The Gauss-Manin connection $\nabla_{\frac{\partial}{\partial t}}=\frac{d}{df}$ on $H$ takes the form
$$
\frac{d}{df}:H\rightarrow H_C,\ \omega\mapsto \omega':=\frac{d\omega}{df},
$$
where $H_C$ is the localization of $H$ on the multiplicative group generated by $t-c_i,\ i=1,2,
\ldots,\mu$ (see \cite{mov}). 

Let $\F=\F(df)$ be the foliation in $U$ with the first integral $f$.  
We consider the holomorphic foliation 
\begin{equation}
\label{df+omega}
\F_\epsilon: df+\epsilon\omega,\ \epsilon\in(\C,0),\ \omega\in \Omega_U^1.
\end{equation}
 Let $b$ be a regular point of $f$, $p\in U_b$ and 
$\Sigma$ be a transversal section at $p$ to $\F(df)$ parameterized by the image $t$ of $f$. 
Let also $\delta\in G:=\pi_1(U_b,p)$ and $h_\epsilon(t):\Sigma\rightarrow \Sigma$ be the holonomy
of $\F_\epsilon$ along the path $\delta$.
 We write the Taylor expansion of $h_\epsilon(t)$ in $\epsilon$:
$$
h_\epsilon(t)-t=M_1(t)\epsilon+M_2(t)\epsilon^2+\cdots +M_i(t)\epsilon^i+\cdots,\  
M_i(t):=\frac{1}{i!}\frac{\partial^ih_\epsilon}{\partial \epsilon^i}\mid_{\epsilon=0}.
$$
$M_i$ is called the $i$-th Melnikov function of the deformation along the path $\delta$.
 For $\delta_1,\delta_2\in G$ we denote by
 $(\delta_1,\delta_2)=\delta_1\delta_2\delta_1^{-1}\delta_2^{-1}$ the commutator
 of $\delta_1$ and $\delta_2$ and for two sets $A,B\subset G$ by $(A,B)$ we mean
 the group generated by $(a,b),\ a\in A, b\in B$. 
 Let
 $$
 G_r:=(G_{r-1},G),\ r=1,2,3,\ldots, \ G_1:=G.
 $$
Using the methods introduced in \cite{ga06} it can be proved that 
if $\delta\in G_k$ then $M_1=M_2=\cdots=M_{k-1}=0$ and  
\begin{equation}
\label{baazdele}
M_k(t)=\int_{\delta_t}\underbrace{\omega_1(\omega_1(\cdots(\omega_1(}_{k-1 \hbox{ times } (}
\omega_1)'\cdots)')'.
\end{equation}
In particular, for $k=2$ and $\delta=(\delta_1,\delta_2)$ we have
\begin{equation}
\label{4mar06}
M_2(t)=\int_{(\delta_1,\delta_2)}\omega\omega'=\det\mat { \int_{\delta_1}\omega} { \int_{\delta_2}\omega}{ \int_{\delta_1}\omega'}{ \int_{\delta_2}\omega'}. 
\end{equation}
\section{Main theorem}
\label{second}
%

For $\omega\in\Omega_U^1$ define $\deg(\omega)$ to be the pole order of $\omega$ along $D$, where $D$ is the compactification divisor of $U$ as it is explained in \S \ref{first}. 
\begin{theo}
\label{commuting}
In the deformation (\ref{df+omega}) with $\deg(\omega)\leq \deg(df)$  assume that 
$f$ is a strongly tame polynomial. Consider $\delta_1,\delta_2\in \pi_1(\Fi_b,p)$ such that
 the corresponding cycles in $H_1(\Fi_b,\Z)$  vanishes along the paths which do not
 intersect each other except at $b$.
 Also assume that 
 \begin{equation}
\label{faca}
\forall \delta\in H_1(\Fi_b,\Z) 
\end{equation}
$$
\langle \delta,\delta_1\rangle =0 \hbox{ or } 
\langle \delta,\delta_2\rangle =0 \hbox{ or }
\exists \delta' \in H_1(\Fi_b,\Z)\hbox{ s.t. }  
\langle \delta_1,\delta\rangle\langle \delta_2,\delta'\rangle-
\langle \delta_2,\delta\rangle\langle \delta_1,\delta'\rangle\not= 0.
$$
If the deformed monodromies along $\delta_1$ and $\delta_2$ commute
then $\F_\epsilon$ has a first integral.
\end{theo}
Note that 
if 
$\langle \delta_1,\delta_2\rangle\not=0$ then the condition in (\ref{faca}) is fulfilled.  
Theorem \ref{12.4} is therefore a special case of Theorem \ref{commuting}. 

\begin{lemm}
\label{11nov05}
Consider a strongly tame morphism $f$, a differential $1$-form $\omega\in \Omega^1_U$ and 
a family of vanishing cycles $\delta=\delta_t$ such that 
$P(t):=\int_\delta \omega$ has the following property: At each $c\in C$, $P$ can be 
written locally in the form $P(t)=(t-c)^\alpha \cdot p(t)$ for some $\alpha\in \C$ and a single-valued holomorphic
function $p$ in $(\C,c)\backslash \{c\}$. 
Then $\omega$ is a relatively exact 1-form and so $\int_{\delta}\omega$ is identically zero.   
\end{lemm}
\begin{proof}
Take a vanishing cycle $\delta'$ with the corresponding critical value $c'$ of $f$ and
the vanishing path $\gamma'$. If $\langle \delta,\delta'\rangle\not =0$ then 
by the Picard-Lefschetz formula along the path $\gamma'$ and for the cycle $\delta$ we have:
\begin{equation}
\label{12.4.06}
\int_{\delta'}\omega=c_{\delta'}P(t),
\end{equation}
where $c_{\delta'}$ is some constant depending on $\delta'$.
Since the Dynkin diagram of $f$ is connected, the equality 
in (\ref{12.4.06}) holds  for 
all vanishing cycles $\delta'$. 
Let $\delta$ be a vanishing cycle along the  path $\gamma$ in the critical value $c$. 
Since $\langle \delta,\delta\rangle=0$, 
the value of the integral $\int_{\delta }\omega$ after the monodromy along the path $\gamma$ 
and around $c$
does not change  and so the corresponding $\alpha$ must be integer.

We conclude that $\int_{\delta}\omega$ for any vanishing cycle $\delta$, is a  single-valued function
 in $\C\backslash
C$. Using the Picard-Lefschetz formula for two vanishing cycles $\delta_i,\delta_j$ with non zero intersection
number, we conclude that $\omega$ is a relatively exact $1$-form.
\end{proof}

{\it Proof of Theorem \ref{commuting}:}
The first Melnikov function associated to the path $(\delta_1,\delta_2)$ is trivially zero.
By  (\ref{4mar06})
we have:
\begin{equation}
\label{elzein}
\frac{\int _{\delta_1} \omega'}{\int_{\delta_1}\omega}=
\frac{\int _{\delta_2}\omega'}{\int_{\delta_2}\omega}.
\end{equation}
If for a continuous family of vanishing cycles $\delta$, we have $\int_{\delta}\omega=0$
then by Proposition \ref{6mar06} the $1$-form $\omega$ is relatively exact and so 
$\omega=Pdf+dQ$ for two meromorphic function in $M$ with poles along $D$. The hypothesis
$\deg(\omega)\leq \deg(df)$ implies that $P=0$ and so $\omega=dQ$.
 
 Therefore, we can assume that that $\int_{\delta_i}\omega,\ i=1,2$ are  not identically zero.
Then the multi-valued function (\ref{elzein}) is well-defined.   
We denote it by $P$ and claim that $P$ is a rational function.
Since integrals have finite growth at critical points and at infinity, it is enough to prove that
$P$ is single-valued. By the hypothesis on the vanishing paths $\gamma_i,\ i=1,2$ of $\delta_i$, 
we can put $\gamma_i$ inside a distinguished system of paths $\Gamma$. 
Let $c\in C$ and $\delta$ be the corresponding 
vanishing cycle along the path $\gamma\in \Gamma$. 
By the Picard-Lefschetz formula 
along the path $\gamma$ 
we have:
$$
\frac{\int _{\delta_1} \omega'+r_1\int_{\delta}\omega'}{\int_{\delta_1}\omega+r_1\int_{\delta}\omega}=
\frac{\int _{\delta_2}\omega'+r_2\int_{\delta}\omega'}{\int_{\delta_2}\omega+r_2\int_{\delta}\omega},
$$
where $r_i:=\langle \delta_i, \delta\rangle, \ i=1,2$.
This and (\ref{elzein}) imply together that either 
$P(t)=\frac{\int _{\delta} \omega'}{\int_{\delta}\omega}$  or $\int_{r_1\delta_2-r_2\delta_1}\omega=0$.
If one of $r_i$'s is zero then $P$ is single-valued along $\gamma$. If both are non-zero then the second case cannot happen because of Proposition \ref{6mar06} and the 
hypothesis in (\ref{faca}). In the first case
we conclude that $P$ is again single-valued in a neighborhood of  $\gamma$.

We have proved that $P$ is a rational function. 
Now $\ln (\int_{\delta_1}\omega)'=P(t)$ and so 
$$
\int_{\delta}\omega=e^{\int P(t)dt}=Q\prod_{c\in K}(t-c)^{\alpha_c},\ Q\in\C(t),\ \alpha_c\in\C
$$
where $K$ is a finite subset of $\C$.
Lemma \ref{11nov05} finishes the proof.
$\Box$

Concerning  the comments after Theorem \ref{12.4} note that for a hyperplane section $D$ 
of $M$ we have the long exact sequence
$$
\cdots \rightarrow H_2(M,\Q)\stackrel{s_1}{\rightarrow} H_2(M,U,\Q)\rightarrow H_1(U,\Q)\rightarrow H_1(M,\Q)\rightarrow \cdots 
$$
and the Leray-Thom-Gysin isomorphism $s_2: H_2(M,U,\Q)\to H_0(D,\Q)\cong \Q$. The composition $s_2\circ s_1$
is the intersection with $D$ and so $s_1$ is not the zero map. This implies that if $H_1(M,\Q)$ is zero
then $H_1(U,\Q)=0$.  


\def\cprime{$'$} \def\cprime{$'$} \def\cprime{$'$}

\begin{tabular}{ll}
Hossein  Movasati  & \qquad Isao Nakai\\
IMPA-Estrada D. Castorina, 110 & \qquad Ochanomizu University\\
Jardim Bot\^anico  & \qquad  Department of Mathematics\\
Rio de Janeiro - RJ  & \qquad   2-1-1 Otsuka, Bunkyo-ku\\
CEP. 22460-320   & \qquad Tokyo 112-8610\\
Brazil &  \qquad Japan\\
{\tt hossein@impa.br}& \qquad {\tt nakai@math.ocha.ac.jp} 
\end{tabular}

\end{document}